\documentclass[11pt]{amsart}
\usepackage{setspace}
\usepackage{graphicx}
\usepackage{amsthm,amsmath,verbatim,amssymb}
\usepackage[colorlinks=true,linkcolor=black,citecolor=black,urlcolor=black]{hyperref}
\usepackage{subcaption}
\usepackage{caption}
\usepackage{color}
\usepackage{arydshln}
\usepackage{mathrsfs}
\usepackage{amscd}
\usepackage{amsfonts}
\usepackage{amstext}
\usepackage[all]{xy}
\usepackage{chemarrow}
\usepackage{indentfirst}
\usepackage{graphicx}
\usepackage{amsfonts}
\usepackage{ulem}

\newtheorem{thm}{Theorem}[section]
\newtheorem{lem}[thm]{Lemma}

\newtheorem{cor}[thm]{Corollary}

\newcommand{\beq}{\begin{eqnarray}}
\newcommand{\eeq}{\end{eqnarray}}
\newcommand{\beqq}{\begin{eqnarray*}}
\newcommand{\eeqq}{\end{eqnarray*}}
\newcommand{\bc}{\begin{center}}
\newcommand{\ec}{\end{center}}

\newcommand{\Tube}{{\rm Tube}}
\def\definedas{\stackrel{\Delta}{=}}

\makeatletter
\@addtoreset{equation}{section}
\makeatother

\begin{document}
\title[Critical radius]{Critical radius and supremum of random spherical harmonics (II)}

\author[ Feng]{Renjie Feng}
\author[ Xu]{Xingcheng Xu}
\author[ Adler]{Robert  J.\ Adler}

\address{Beijing International Center for Mathematical Research, Peking University, Beijing, China}
\email{renjie@math.pku.edu.cn}

\address{School of Mathematical Sciences, Peking University,
Beijing, China}
\email{xuxingcheng@pku.edu.cn}
\address{Andrew and Erna Viterbi Faculty of Electrical Engineering, Technion-Israel Institute of Technology, Haifa 32000, Israel}
\thanks{Research of RJA supported in part by URSAT, ERC  Advanced Grant 320422.}
\email{radler@technion.ac.il}

\date{\today}
\maketitle


\begin{abstract}
We continue the study, begun in \cite{FA}, of the critical radius of embeddings, via deterministic spherical harmonics,  of  fixed dimensional  spheres into higher dimensional ones, along with the associated problem of the distribution of the suprema of  random spherical harmonics.
Whereas  \cite{FA} concentrated on spherical harmonics of a common degree, here we extend the results to mixed degrees, obtaining larger lower bounds on critical radii than we found previously.
\end{abstract}

\section{Introduction}
The spherical harmonics, of level $n \geq 1$, on the $d$-dimensional unit sphere $S^d$, are the collection of
 eigenfunctions $\{\phi_j^{n,d}\}_{j=1}^{k_n^d}$ of the Laplacian $\Delta_{g_{S^d}}$ on $S^d$, satisfying
\begin{equation}
\Delta_{g_{S^d}} \phi^{n,d}_j(x) \ = \-n(n+d-1)\phi_j^{n,d}(x),
\end{equation}
where $k_n^d$ is defined below at  \eqref{eq:kld}.

In \cite{FA}, we studied  the map
 \begin{equation}
 \label{af}
 \tilde i_n^d:\,\,\, S^d\to S^{k^d_n-1},\ \quad x\to \sqrt{\frac{s_d}{k^d_n}}\left(\phi_{1}^{n,d}(x), \cdots  ,\phi_{k_n^d}^{n,d}(x)\right),
 \end{equation}
defined by the spherical harmonics of level $n$, where, $s_d$ denotes the surface area of the unit sphere $S^d$.
For large enough $n$,  the image   $\tilde i^d_n(S^d)$ is diffeomorphic  to $S^d$ if $n$ is odd, and  to $\mathbb RP^d$ if $n$ is even.

As shown in \cite{FA}, and explained again below, the behavior of this {\it deterministic} map has significant implications for {\it random} spherical harmonics of the form
\beq
\tilde \Phi^d_n(x)\ = \ \sum_{j=1}^{k^d_n} a_j^n\phi_j^{n,d}(x),
\label{eq:rob1}
\eeq
where the $a^n_j$ are either standard Gaussian variables (in which case we talk about  the ``Gaussian ensemble") or $(a^n_1,\dots, a^n_{k^d_n})$ is uniform on $S^{k^n_d-1}$ (in which case we talk about the ``spherical ensemble"). In particular, one of the most important aspects of the deterministic mapping vis a vis the random process is the critical radius, or reach, of its image in the ambient sphere. (For a definition of the notion of critical radius, which plays an important role in  integral geometry and Weyl's tube formula, and for some of its basic properties relevant to our setting, see Section 3 of \cite{FA}.)

What made the  study in \cite{FA} most interesting is the fact that the pull back of the standard round metric under $\tilde i_n^d$    grows with order  $n^2$,  indicating that the image $\tilde i^d_n(S^{d})$ in $S^{k^d_n-1}$  becomes  more and more `twisted' as $n$ grows.
Intuitively, if  a sequence of sets become more twisted  in their ambient space, it seems natural that their critical radii, as a measure of smoothness, will tend to zero. (Think of the critical radii of the graphs of $(x, \sin nx)$ in $\mathbb R^2$, which  tend to 0 as $n\to\infty$.)  Rather surprisingly,  the main results of \cite{FA} showed that there is a lower bound for the critical radius of  the $\tilde i^d_n(S^{d})$ in $S^{k^d_n-1}$, as $n\to\infty$. In some sense, this was  a result of  competition between the `twistiness' of the image and the  `extra space' available as the ambient spaces $S^{k^d_n-1}$ changed.

As a direct consequence of these deterministic results \cite{FA} derived an explicit formula for the distribution of the suprema of the  random spherical harmonics
$\tilde \Phi^d_n$ of \eqref{eq:rob1}
under the spherical ensemble, by exploiting Weyl's tube formula. (For the Gaussian ensemble,   see \ \cite{AT,S1, TK} on the connections
between spherical and Gaussian ensembles, Weyl's tube formula, suprema of random fields, and the expected Euler characteristic of excursion sets.)

The aim of the present paper is to extend the analysis of \cite{FA} to a related, but somewhat different embedding,  given by
 the deterministic map
\begin{equation}\label{id111}
i_n^d:\ S^d\to S^{\pi^d_n-1},\qquad  x\to  \sqrt{\frac{s_d}{\pi^d_n}}\left(\phi_{j}^{\ell,d}(x)\right)_{\ell=0,\dots,n, \ j=1,\dots,k_\ell^d},
\end{equation}
where $\pi^d_n=\sum_{\ell=0}^nk_{\ell}^d$. For large enough $n$ (either  odd or even) this map is an embedding; viz.\ $i_n^d( S^d)\cong S^d$.
(cf.\ \cite{P,Ze3}.)
Following the ideas and proofs in \cite{FA}, we  will prove the existence of a lower bound for the critical radius  of $i_n^d(S^{d})$ in $S^{\pi^d_n-1}$. This will  allow us to also  derive an exact  formula for the distribution  of the suprema of  the  family of spherical ensemble, random spherical harmonics,
\beq
\Phi^d_n(x)\ = \ \sum_{\ell=0}^n\sum_{j=1}^{k^d_\ell} a_j^\ell\phi_j^{\ell,d}(x),
\label{eq:rob11}
\eeq
with $(a_j^\ell)_{\ell,j}$ uniform on $S^{\pi^d_n-1}$. (cf.\ Section \ref{sec:22} below  for details.)

The differences between the two random processes $\Phi^d_n$ and
$\widetilde\Phi^d_n$ are subtle but important, and best understood in spectral terms. For a fixed $\ell$, all the spherical harmonics
$\phi^{\ell,d}_j$ are associated with the same eigenvalue, often called `frequency'. In these  terms, $\widetilde\Phi^d_n$ is a single, or `pure' frequency random process, whereas the spectrum of $\Phi^d_n$ contains (a discrete collection of) frequencies between 0 and
$n(n+d-1)$. In terms of the original motivation for studying random spherical harmonics (the `Berry conjecture' of \cite{Berry}), mixed spectra processes play a more central role than  pure spectra ones. This is one of the main motivations behind the current paper.

The second motivation is a purely mathematical one.  As described above, the existence of a lower bound for the critical radius of
$\widetilde i_n^d(S^d)$ involves a delicate balance between the contortions of the set itself and the high dimension of the space in which it is embedded. It was not at all clear,  a priori, how these results would extend to the reach of  $i_n^d(S^d)$. On the one hand, the addition of lower frequencies might be expected to have a smoothing effect, so that, together with the larger dimension of the ambient space, ($\pi_n^d$ rather than $k_n^d$) it seemed reasonable to expect a larger critical radius for $i_n^d(S^d)$ than for
$\widetilde i_n^d(S^d)$. On the other hand, it would also be not unreasonable that the highest frequencies dominate, and that there be no difference between the two cases. Back-of-the-envelope calculations are not accurate enough to enable one to  make even a reasonable conjecture as to which is the case, but the detailed calculations of this paper show that the first scenario is, in fact, the correct one.

In the following section we  set up some necessary  notations, and then formally state all the main results of the paper. Section
\ref{sec1} provides proofs for the case $d=2$, and Section \ref{eeee} describes the extensions needed for the general $d$ cases.

\section{Main results}\label{sec:main}

\subsection{Spherical harmonics and the deterministic embedding}
Let the unit sphere $S^d$ be equipped with the round metric $g_{S^d}$, and write  $\Delta_{g_{S^d}}$ for the associated Laplacian.
The spherical harmonics  $\phi^{\ell,d}_j$ of level $\ell\geq 0$ are  the eigenfunctions
\begin{equation}
\Delta_{g_{S^d}} \phi^{\ell,d}_j(x) \ = \-\ell(\ell+d-1)\phi_j^{\ell,d}(x).
\end{equation}
Write  $\mathscr{H}^d_\ell$ for the eigenspace spanned by these eigenfunctions. Then  the dimension of  $\mathscr{H}^d_\ell$  is
\beq
k^d_\ell\ \definedas \ \frac{2\ell+d-1}{\ell+d-1}{\ell+d-1\choose d-1}.
\label{eq:kld}
\eeq
Since  $L^2(S^d)=\oplus_{\ell\geq 0}\mathscr{H}_\ell^d$,
if we normalize the eigenfunctions so that their $L^2$-norm  is  $1$, the expansion of $L^2(S^d)$ functions in the  orthonormal basis of  spherical harmonics provides a natural generalization of  Fourier series expansions.

Now set
$$\mathscr{H}^d_{\ell\leq n}\ \definedas\ \bigoplus_{\ell=0}^n\mathscr{H}_\ell^d$$
to denote the space of  spherical harmonics of degree at most $n$. Then the dimension $\pi^d_n$ of $\mathscr{H}^d_{\ell \leq n}$ is \cite{AH} \begin{equation}\label{dimn}
\pi^d_n\  \definedas \ \sum_{\ell=0}^nk_{\ell}^d\ =\ \frac{2n+d}{d}{n+d-1\choose d-1} .
\end{equation}

Consider now the map
\begin{equation}\label{id111}
i_n^d:\ S^d\to \mathbb R^{\pi^d_n},\ x\to  \sqrt{\frac{s_d}{\pi^d_n}}\left(\phi_{j}^{\ell,d}(x)\right)_{\ell=0,1,\cdots,n, \  j=1,\cdots,k_\ell^d}.
\end{equation}
For $n$ large enough, this map is an embedding \cite{P,Ze3}  i.e., $i_n^d( S^d)\cong S^d$. Furthermore, it follows from the
 properties of  spectral projection kernels that the $\mathbb R^{\pi_n^d}$ norm of $i_n^d(x)$ is identically $1$, so that $i_n^d$ is actually a map between spheres; viz.
\begin{equation}\label{id2}
i_n^d:\,\,\, S^d\to  S^{\pi_n^d-1}.
\end{equation}
In addition, the pull-back of the Euclidean metric satisfies  \cite{Ze3} \begin{equation}(i_n^d)^*(g_E)\cong c_dn^2 g_{S^d}\end{equation}
where $c_d$ is a constant  depending  only on $d$.

Our interest in this section lies in the critical radius of $i_n(S^d)$ in $\mathbb R^{\pi_n^d}$.

 Recall that if
 $M$ is a smooth manifold  embedded in an ambient manifold $\widetilde M$, then
the local critical radius, or reach, at a point $x\in M$ is the furthest distance one can travel, along any geodesic in $\widetilde M$  based at $x$ but normal to $M$ in $\widetilde M$, without meeting a similar vector originating at another point in $M$. The (global) critical radius of $M$ is then the infimum of all the local ones.
We refer the reader to Section 3 of \cite{FA}
for additional background and for formal definitions.

We can now state our first result.

%
%
%
%
%

\begin{thm}\label{1}

For  sufficiently large  $n$, the critical radius  of the embedding $i_n(S^d)$ in $\mathbb R^{\pi_n^d}$ has a strictly positive, uniform in $n$, lower bound which  depends only on $d$.
\end{thm}
Let
$$
\Tube_{\pi_n^d}(i_n^d(S), \rho) \ \definedas \ \left\{
x\in\mathbb R^{\pi_n^d} :\ \min_{p\in i_n^d(S)} \|x-p\|\, \leq \, \rho
\right\}
$$
 be the tube around $i_n(S^d)$ in $\mathbb R^{\pi_n^d}$, where $\rho$ is  less than the critical radius of
 $i_n^d(S)$  in $\mathbb R^{\pi_n^d}$. Then,  by \eqref{id2}, the intersection
 $
 \Tube_{{\pi_n^d}}(i_n^d(S), \rho)\cap  S^{\pi_n^d-1}
 $
 will be a tube of  $i_n(S^d)$ in  $S^{\pi_n^d-1}$ without self-intersection. This fact immediately implies

\begin{cor}\label{dlll}
	Theorem \ref{1} continues to  hold, with a similar lower bound, when  $i_n(S^d)$ is considered as an embedding in $S^{\pi_n^d-1}$.
	 \end{cor}

\subsection{Random spherical harmonics and exceedence probabilities}
\label{sec:22}
In this section we turn our attention to the random spherical harmonics under the spherical ensemble. These are defined by
\beq
\Phi^d_n(x) \definedas \  \sum_{\ell=0}^n\sum_{j=1}^{k^d_\ell} a_j^\ell\phi_j^{\ell,d}(x),\label{eq:rob2}
\eeq
where the random vector $a=(a_j^\ell)_{\ell=0,\dots,n,\ j=1,\dots, k_\ell^d}$ is distributed uniformly over the sphere $S^{\pi_n^d-1}$.
%
%
%
%
As opposed to the simpler
 $\widetilde \Phi^d_n$ of \eqref{eq:rob1}, $ \Phi^d_n$  has a broad spectrum.


 If we now let $\rho_d$ denote the uniform lower bound for the critical radius of $i_n^d(S^d)$ in the ambient space $S^{\pi_n^d-1}$ appearing in Corollary  \ref{dlll}, then this corollary and
 the same arguments as adopted in Section 6 of \cite{FA} prove the following result.

\begin{thm}\label{34}Let $\Phi_n^d$ be the random spherical harmonics under the spherical ensemble as defined above.
Then there exists constants $\rho_d>0$ such that, for sufficiently large $n$, and for all $u>\sqrt{\frac{\pi_n^d}{s_d}}\cos(\rho_d)$,
\begin{equation}
\begin{split}
& \mathbb P_{\mu^d_n}\left\{\sup_{S^d}\Phi^{d}_n(x)>u\right\}\\
&\ \ \ \  \ = \ \frac 1{s_{\pi_n^d-1}}\sum_{j=0}^d f_{\pi_n^d,j}\left(\cos^{-1}\left(u/\sqrt{\frac{\pi_n^d}{s_d}}\right)\right) \left(\frac{n(n+d)}{d+2}\right)^{j/2}\mathcal L_j(S^d),
\end{split}
\end{equation}
where $s_{\pi_n^d-1}$ is the surface area of the unit sphere $S^{\pi_n^d-1}$,  the $f_{\pi_n^d, j}(\rho)$ are   explicit functions given in Theorem
10.5.7 of \cite{AT} and {\rm (6.6)} of \cite{FA},  and the $\mathcal L_{j}(S^d)$ are the standard $j$-th Lipschitz-Killing curvatures  of the unit sphere $S^d$, given explicitly, for example,  in {\rm (6.10)} of \cite{FA}.
\end{thm}

Note that although Theorem \ref{34} gives an exact result under the spherical ensemble, analogous (but approximate) results can also
be formulated under the Gaussian ensemble (i.e.\ when the $a_j^\ell$ are all independent, standard, Gaussian random variables). This follows from the rich literature   relating mean Euler characteristics of excursion sets and exceedence probabilities  for Gaussian processes; e.g.\  \cite{AT, CX, S1, TK, VM}. We will not repeat this here.

\subsection{Remark}

Before turning to proofs, we want to make a comment about the possibility of extending the results above to a setting of  manifolds other than the sphere.

Let $(M, g)$ be an $m$-dimensional Riemann manifold and $\Delta_g$ its Laplace-Beltrami operator.
Fix $\lambda > 0$ (typically large) and consider the map $$i_\lambda: M\to \mathbb R^{d_\lambda},\,\,\,\, x\to \frac 1{\Pi_\lambda(x,x)} (\phi_1(x),\cdots, \phi_{d_\lambda}(x)),$$
where the $\phi_j$, $j\leq d_\lambda$,  are the eigenfunctions of $\Delta_g$ corresponding to  eigenvalues  less than $\lambda$,  and $\Pi_\lambda$ is the spectral projection kernel.

Following the strategy in the proofs to follow, in order to establish a uniform estimate for the critical radius of $i_\lambda(M)$ in  $\mathbb R^{d_\lambda}$, one needs to consider the on and off-diagonal estimates of the spectral projection kernel  $\Pi_\lambda$ and its derivatives.
In the case of the round sphere, we have  explicit expressions for the kernel by the Christoffel-Darboux formula. Furthermore, the local behavior of the kernel is encoded by a Hilb's type formula, giving Mehler-Heine estimates.
On the other hand, for general manifolds, globally, we only have H\"ormander's classical estimate   (Theorem 21.1 in \cite{SHu}) that
$$
|\Pi_\lambda(x,y)|\leq c(1+\sqrt  \lambda)^{m-1},
$$
 if the pair $(x,y)$ belongs to a compact set in $M \times M$ disjoint from the diagonal. On the diagonal, there is a universal  limit for the spectral projection kernel in the ball of the length scale $1/\sqrt \lambda$; viz.
 $$
 \frac 1{\lambda^{m/2}}\Pi_\lambda\left(p+\frac x {\sqrt \lambda}, p+\frac y {\sqrt \lambda}\right)
 $$
  has a $\lambda\to\infty$ limit, which can be expressed in term of  Bessel functions \cite{Sz}.

  However, there are regions for which we do not have general bounds on $\Pi_\lambda(x,y)$, a typical example being  when $\|x-y\|$ lies in the length scale $[{\lambda^{-1/3}},1]$.  It is this lack of useful bounds that makes tackling the same problems, but on general manifolds, currently intractible.  Nevertheless,  we hope to establish extensions of this kind in the future.

\section{Proof of Theorem \ref{1} for $S^2$}\label{sec1}
In this section, we will first collect some standard results regarding  spherical harmonics on the $2$-sphere that we will need later, and then prove Theorem \ref{1} for the $2$-sphere.

\subsection{Spectral projection kernels}\label{sec11}

We will drop the index 2 in this section whenever it does not lead to ambiguities.
Thus $\mathscr{H}_{\ell\leq n}$ is now the space of spherical harmonics of $S^2$ of degree at most $n$. Its dimension is $\pi_n=(n+1)^2$ by \eqref{dimn}.  The spectral projection kernel  is now given, for $x,y\in S^2$, by the  Christoffel-Darboux formula \cite{AH},
$$
K_n(x,y)\ = \ \sum_{\ell=0}^n\sum_{j=-\ell}^\ell \phi_j^\ell(x)\phi_j^\ell(y)\ =\ \frac{n+1}{4\pi}P_n^{(1,0)}(\cos\Theta(x,y)),
$$
where $\Theta(x,y)$ is the angle between the vectors $x$ and $y$ on $S^2$, and
 $P_n^{(1,0)}$ is a Jacobi polynomial. In general, the Jacobi polynomials are defined by
$$P_n^{(\alpha,\beta)}(x)\ = \ \sum_{s=0}^n \binom{n+\alpha}{s}\binom{n+\beta}{n-s}\left(\frac{x-1}{2}\right)^{n-s}
\left(\frac{x+1}{2}\right)^{s}.$$
We will need the following facts about $P^{(1,0)}$:
$$P_n^{(1,0)}(1) \ = \ n+1,\qquad  P_n^{(1,0)}(1)'\ = \ \frac{n(n+1)(n+2)}{4}.$$
We will also need the fact that, on the diagonal, the kernel $K_n$ is given by
$$K_n(x,x)\ =\ \frac{n+1}{4\pi}P_n^{(1,0)}(1)\ =\ \frac{(n+1)^2}{4\pi}.$$
Defining the normalized kernel
$$\Pi_n(x,y)\ \definedas \ \frac {4\pi}{(n+1)^2}K_{n}(x,y)\ = \ \frac {1}{n+1}P_n^{(1,0)}(\cos\Theta(x,y)),$$
we have that the norm of   $i_n(x)$, defined at  \eqref{id111}, is given by
$$\|i_n(x)\|^2\ =\ \frac {4\pi}{(n+1)^2}\sum_{j,\ell}|\phi_{j}^\ell(x)|^2\ =\ \Pi_n(x,x)\ =\ 1,$$
so that $i_n$ is actually a map
\begin{equation}\label{iii}i_n:\ S^2 \to S^{n^2+2n}.\end{equation}
Following the  computations in \cite{FA}, the pull back of the Euclidean metric under this mapping is
\begin{equation}\label{metric}
g_n\ = \ i_n^*(g_E) \ = \ \frac{n^2+2n}{4} g_{S^2}.
\end{equation}

\subsection{Critical radius of $i_n(S^2)$}\label{ddddd}
A useful formula  for the critical radius of a smooth manifold embedded in Euclidean space was  derived in \cite{TK}. For the case of the embedding of $i_n(S^2)$ this gives the critical radius as


\begin{equation}\label{critical}
r_{n} \ =\ \inf_{x,y\in S^2} \frac{\|i_n(x)-i_n(y)\|^2}{2 \|P^{\perp}_{i_n(y)}(i_n(x)-i_n(y))\|},
\end{equation}
where $P^{\perp}_y(x-y)$ is the projection of the $x-y$ to the normal bundle at $y$.


Following the  calculations in Section 3 of \cite{FA},  this can be rewritten as
\begin{equation}r_{n}\ =\ \inf_{\theta\in [0,\pi]} \frac{1-\frac{1}{n+1}P_n^{(1,0)}(\cos\theta)}{\sqrt{2-\frac{2}{n+1}P_n^{(1,0)}(\cos\theta)-\frac{1}{n+1}\frac{[ (P_n^{(1,0)}(\cos\theta))'\sin\theta]^2}{(P_n^{(1,0)}(1))'}}}.
\label{eq:rob5}
\end{equation}

\subsection{Proof of Theorem \ref{1} for $S^2$}\label{dddd}
For simplicity, we use  $P_n$ to denote the Jacobi polynomial $P_n^{(1,0)}$  for the remainder of this section. The  proof is then based on   classical asymptotic estimates for   Jacobi polynomials.

We start with  an asymptotic formula of Hilb's type (\cite{Sz}, Theorem 8.21.12):
\begin{equation}\label{Hilb}
P_n(\cos \theta)
\ = \ \left(\sin \frac{\theta}{2}\right)^{-1}\left\{\left(\frac{\theta}{\sin \theta}\right)^{\frac{1}{2}}J_1((n+1)\theta)+R_{1,n}(\theta)\right\},
\end{equation}
where
\begin{equation}
      R_{1,n}(\theta)\ = \
      \begin{cases}
        \theta^{3}O(n),\ \ \ \ \ \  0\leq \theta\leq c/n,\\
        \theta^{\frac{1}{2}}O(n^{-\frac{3}{2}}),\ \  c/n\leq \theta\leq \pi-\epsilon,\\
      \end{cases}
\end{equation}
and $c,$ $\epsilon$ are uniform constants, independent of $n$.

We also have the Darboux formula (\cite{Sz}, Theorem 8.21.13):

\begin{equation}\label{Darboux}
P_n(\cos\theta)\ = \ \frac{1}{\sqrt{n}}k(\theta)\cos\left((n+1)\theta-\frac{3}{4}\pi\right)+R_{2,n}(\theta),
\end{equation}
if $c'/n\leq \theta \leq \pi-c'/n$, where $c'$ is a uniform constant and
\beq
k(\theta)&=&\frac{1}{\sqrt{\pi}}\left(\sin\frac{\theta}{2}\right)^{-3/2}\left(\cos\frac{\theta}{2}\right)^{-1/2},\\
R_{2,n}(\theta) &=& \frac{1}{\sqrt{n}}k(\theta)\frac{1}{n\sin\theta}O(1).
\eeq

Near the right end point $\pi$, the asymptotic behavior of the Jacobi polynomials is given by the Mehler-Heine formula
\begin{equation}\label{MH}
\lim_{n\to\infty}P_n\left(\cos\left(\pi-\frac{x}{n}\right)\right)\ = \ J_0(x),
\end{equation}
where the limit is uniform on compact subsets of $\mathbb{R}$ (see equation (3) in  \cite{bmo}) and the $J_\nu$ are the Bessel functions of the first kind of order $\nu$.

In order to prove Theorem \ref{1}, we divide $[0, \pi]$ into the four subintervals
\beq
\label{eq:ints}
[0,c/n],\ \  [c/n, n^{-4/5}],\ \  [n^{-4/5}, \pi-c'/n],\ \  [\pi-c'/n,\pi],
\eeq
and replace  the global  infimum in \eqref{eq:rob5} as a minimium, in self-evident notation, as
\begin{equation}\label{inf4}
\min\left\{I_n, II_n, III_n,IV_n\right\} \ \definedas\  \min\left\{\inf_{[0,c/n]}, \inf_{[c/n,n^{-4/5}]}, \inf_{[n^{-4/5},\pi-c'/n]}, \inf_{[\pi-c'/n,\pi]}\right\} .
\end{equation}
We now treat each of the four infima in \eqref{inf4}, starting  with
$$I_n\ = \ \inf_{[0, c/n]}\frac{1-\frac{1}{n+1}P_n(\cos \theta)}{\sqrt{2-\frac{2}{n+1}P_n(\cos\theta)-\frac{1}{n+1}\frac{ [P'_n(\cos\theta)\sin\theta]^2}{P_n'(1) }}}. $$

To treat this term, we study a rescaling limit via a new parameter $x$, where $x=n\theta$, so that $x\in [0, c]$. By the Hilb's type asymptotic (\ref{Hilb}) on $0\leq \theta\leq c/n$,
\begin{equation*}
         \begin{split}
            \frac{1}{n}P_n\left(\cos \frac{x}{n}\right)
            &=\frac{1}{n}\left(\sin \frac{x}{2n}\right)^{-1}\left(\frac{x/n}{\sin x/n}\right)^{\frac{1}{2}}J_1\left(x+\frac{x}{n}\right)+O(n^{-2}) \\
            &=\frac{1}{n}\left(\frac{x}{2n}+O(n^{-3})\right)^{-1}
            \left(1+O(n^{-2})\right)^{1/2}\\
            &\ \ \ \ \  \ \ \  \times\left(J_1(x)+O(n^{-1})\right)+O(n^{-2}) \\
            &=\frac{2}{x}J_1(x)+O(n^{-1}).\\
         \end{split}
\end{equation*}
Next, for the rescaling of $P'_n\left(\cos\frac{x}{n}\right)$, note the following two standard facts about Jacobi polynomials and Bessel functions:
\beqq
P'_n(x)&=&\frac{n+2}{2}P_{n-1}^{(2,1)}(x),\quad\ \  \text{(\cite{Sz}, (4.5.5))}
\\
\left(\frac{2}{x}J_1(x)\right)' &=&-\frac{2}{x}J_2(x).\qquad \qquad\text{(\cite{Sz}, (1.71.5))}
\eeqq
Applying these facts and the Hilb's asymptotic for $P_{n}^{(2,1)}(\cos\theta)$ given below in  \eqref{jacobi1}, we have
\begin{equation*}
         \begin{split}
            &\frac{1}{n^2}P'_n\left(\cos\theta\right)\left(\sin\theta\right)\\
            &\qquad\ = \ \frac{1}{n^2}\left\{\frac{1}{2}\sin\theta\cdot (n+2)P_{n-1}^{(2,1)}(\cos\theta)\right\} \\
            &\qquad\ = \ \frac{n+2}{n^2}\left(\sin\frac{\theta}{2}\right)^{-1}\left\{\frac{n}{n+1}
            \left(\frac{\theta}{\sin \theta}\right)^{\frac{1}{2}}J_2((n+1)\theta)+O(n^{-2})\right\} \\
            &\qquad\ = \ \frac{n+2}{n^2}\left(\sin\frac{x}{2n}\right)^{-1}\left\{\frac{n}{n+1}
            \left(\frac{x/n}{\sin x/n}\right)^{\frac{1}{2}}J_2\left(x+\frac{x}{n}\right)+O(n^{-2})\right\} \\
            &\qquad\ = \ \frac{2}{x}J_2(x)+O(n^{-1}) \\
            &\qquad\ = -\left(\frac{2}{x}J_1(x)\right)'+O(n^{-1}).
         \end{split}
\end{equation*}
We rescale  $$\frac{1}{n+1}\frac{ [P'_n(\cos\theta)\sin\theta]^2}{P_n'(1)}$$ to obtain
$$
\frac{4\left(\left(\frac{2}{x}J_1(x)\right)'n^2+O(n)\right)^2}{n(n+1)^2(n+2)}
\ =\ 4\left(\left(\frac{2}{x}J_1(x)\right)'+O(n^{-1})\right)^2 .
$$
Hence, as $n\to\infty$, $I_n$ is  asymptotic to
\begin{equation}\label{In}
I_\infty\ =\ \inf_{[0,c]} \frac{1-\frac{2}{x}J_1(x)}{\sqrt{2-\frac{4}{x}J_1(x)-4((\frac{2}{x}J_1(x))')^2 }}.
\end{equation}

As an aside, note that the ratio here is well defined as $x\to0$, as follows from the expansion
$J_1(x)=\frac{1}{2}x-\frac{1}{16}x^3+\frac{1}{384}x^5+O(x^7)$ of the Bessel function around $x=0$.

For the second term, $II_n$, we again use the Hilb's type asymptotic formula \eqref{Hilb}, this time on the interval  $[c/n,n^{-4/5}]\subset [c/n,\pi-\epsilon]$. We also rescale to $x=n\theta$, so that $x\in [c,n^{1/5}]$, and will need the following three  basic properties of Bessel functions (\cite{Sz}, Pages 15--16), for $\nu$ real but  $\nu\neq -1,-2,-3,\dots$.
\beqq
J_\nu(x)&=&\left(\frac{2}{\pi x}\right)^{\frac{1}{2}}\cos(x-\frac{\pi \nu}{2}-\frac{\pi}{4})+O(x^{-\frac{3}{2}}),\quad \text{as}\ x\to \infty,\\
J_\nu(x)&\sim& x^{\nu},\quad \text{as}\ x\to 0^+,\\
J'_1(x)&=&\frac{1}{x}J_1(x)-J_2(x).
\eeqq
Applying these properties, we have, uniformly in $x\in [c,n^{1/5}]$,
\begin{equation}
J_1(x+\frac{x}{n})\ =\ J_1(x)+\frac{x}{n}J'_1(x)+\dots\ =\ J_1(x)+O(n^{-4/5}).
\end{equation}
Consequently,
\begin{equation*}
         \begin{split}
            \frac{1}{n}P_n\left(\cos \frac{x}{n}\right)
            &\qquad\ = \ \frac{1}{n}\left(\sin \frac{x}{2n}\right)^{-1}\left(\frac{x/n}{\sin x/n}\right)^{\frac{1}{2}}J_1(x+\frac{x}{n})+O(n^{-\frac{19}{10}}) \\
            &\qquad\ = \ \frac{1}{n}\left(\frac{x}{2n}+O(n^{-12/5})\right)^{-1}
            \left(1+O(n^{-8/5})\right)^{1/2}\\
            &\qquad\quad \qquad \times\left(J_1(x)+O(n^{-4/5})\right)+O(n^{-\frac{19}{10}}) \\
            &\qquad\ = \ \frac{2}{x}J_1(x)+O(n^{-4/5}).\\
         \end{split}
\end{equation*}
Hence, we have the rescaling limit
\begin{equation}
\frac{1}{n+1}P_n\left(\cos\frac{x}{n}\right)\ =\ \frac{2}{x}J_1(x)+O(n^{-4/5}),
\end{equation}
for $x\in [c,n^{1/5}]$. Similarly, the rescaling of $$\frac{1}{n+1}\frac{ [P'_n(\cos\theta)\sin\theta]^2}{P_n'(1)}$$
will be dominated by the leading term $4((\frac{2}{x}J_1(x))')^2$ for $n$ large enough. Thus  $II_n$ will converge, as $n\to\infty$, to the same expression that we had for $I_n$, viz.
\begin{equation}\label{IIn}
II_\infty\ =\ \inf_{[c,\infty]} \frac{1-\frac{2}{x}J_1(x)}{\sqrt{2-\frac{4}{x}J_1(x)-4((\frac{2}{x}J_1(x))')^2 }}.
\end{equation}

For $III_n$, since $\theta\in [n^{-4/5},\pi-c'/n]\subset [c'/n,\pi-c'/n]$, we can apply the Darboux formula (\ref{Darboux}) on $[n^{-4/5},\pi-c'/n]$ and define
\begin{equation*}
M_n(\theta)\ = \ \frac{1}{\sqrt{n}}
k(\theta)\cos\left((n+1)\theta-\frac{3}{4}\pi\right),
\end{equation*}
giving
\begin{equation*}
P_n(\cos\theta)\ =\ M_n(\theta)+R_{2,n}(\theta).
\end{equation*}
Since $k'(\frac{2\pi}{3})=0$ and $k''(\theta)>0$ for $\theta\in (0,\pi)$, it is simple to check that
\begin{equation*}
\max k(\theta)\ =\ \max\{k(n^{-4/5}),\  k(\pi-c'/n)\}\ =\ O(n^{6/5}),
\end{equation*}
so that
\begin{equation*}
\frac{1}{n}M_n(\theta)\ =\  n^{-3/2}\cdot O(n^{6/5})\ = \ O(n^{-3/10}).
\end{equation*}
Similarly, we have
\begin{equation*}
\frac{1}{n}R_{2,n}(\theta)\ = \ n^{-5/2}\cdot O(n^{2})\ = \ O(n^{-1/2}),
\end{equation*}
so that
$$\frac{1}{n+1}P_n(\cos \theta)\ = O(n^{-3/10}).$$
Note now the fact that, for $c'/n\leq \theta\leq \pi-c'/n$, it follows from bounds on the  derivatives of Jacobi polynomials (\cite{Sz}, P. 236, 8.8.1) that
$$
\frac d{d\theta}{P_n(\cos \theta)}\  =\ n^{\frac{1}{2}}k(\theta)\left\{-\sin((n+1)\theta-3\pi /4)+(n\sin \theta)^{-1}O(1)\right\}.
$$
Also, on the interval $[n^{-4/5},\pi-c'/n]$, we have
\begin{equation*}
\frac d{d\theta}{P_n(\cos \theta)}\ =\ n^{\frac{1}{2}}\cdot O(n^{6/5})\ = \ O(n^{17/10}).
\end{equation*}
This implies that $$\frac{1}{n+1}\frac{ [P'_n(\cos\theta)\sin\theta]^2}{P_n'(1)}\ =\ O(n^{-3/5}).$$
Thus the  $n\to \infty$ limit of $III_n$ is
\begin{equation}\label{IIIn}
III_\infty\ = \ \frac{1}{\sqrt{2}}.
\end{equation}

We treat the final term, $IV_n$, a little differently,  bounding it from below.  Firstly, we have
\begin{equation*}
         \begin{split}
            IV_n
            &=\ \inf_{[\pi-c'/n,\pi]}\frac{1-\frac{1}{n+1}P_n(\cos \theta)}{\sqrt{2-\frac{2}{n+1}P_n(\cos\theta)-\frac{1}{n+1}\frac{ [P'_n(\cos\theta)\sin\theta]^2}{P_n'(1) }}}\\
            &\geq\  \inf_{[\pi-c'/n,\pi]}\sqrt{\frac{1}{2}\left(1-\frac{1}{n+1}P_n(\cos \theta)\right)}
            \\
           &=\  \inf_{[0,c']}\sqrt{\frac{1}{2}\left(1-\frac{1}{n+1}P_n(\cos(\pi-\frac{x}{n}))\right)}.
             \end{split}
\end{equation*}
By the Mehler-Heine formula \eqref{MH}, we have the uniform estimate
$$\ \frac{1}{n+1}P_n(\cos(\pi-\frac{x}{n}))\ \to\  0$$
for $x\in[0,c']$.
Hence, as $n\to\infty$, we have \begin{equation}\label{IVn}
IV_\infty\ \geq\ \frac{1}{\sqrt{2}}.
\end{equation}

The limits and lower bound for $I_n$--$IV_n$ established   above show that the critical radius of   the $i_n(S^2)$, as $n\to\infty$,   have
 a non-zero lower bound in the ambient space $\mathbb R^{(n+1)^2}$. That a similar lower bound holds  in the ambient space $S^{n^2+2n}$ is implied by the discussion before Corollary \ref{dlll}, and this completes  the proof of Theorem \ref{1} for $S^2$.
 \qed

\begin{figure}[ht]
\begin{center}
\includegraphics[width=11cm, height=5cm]{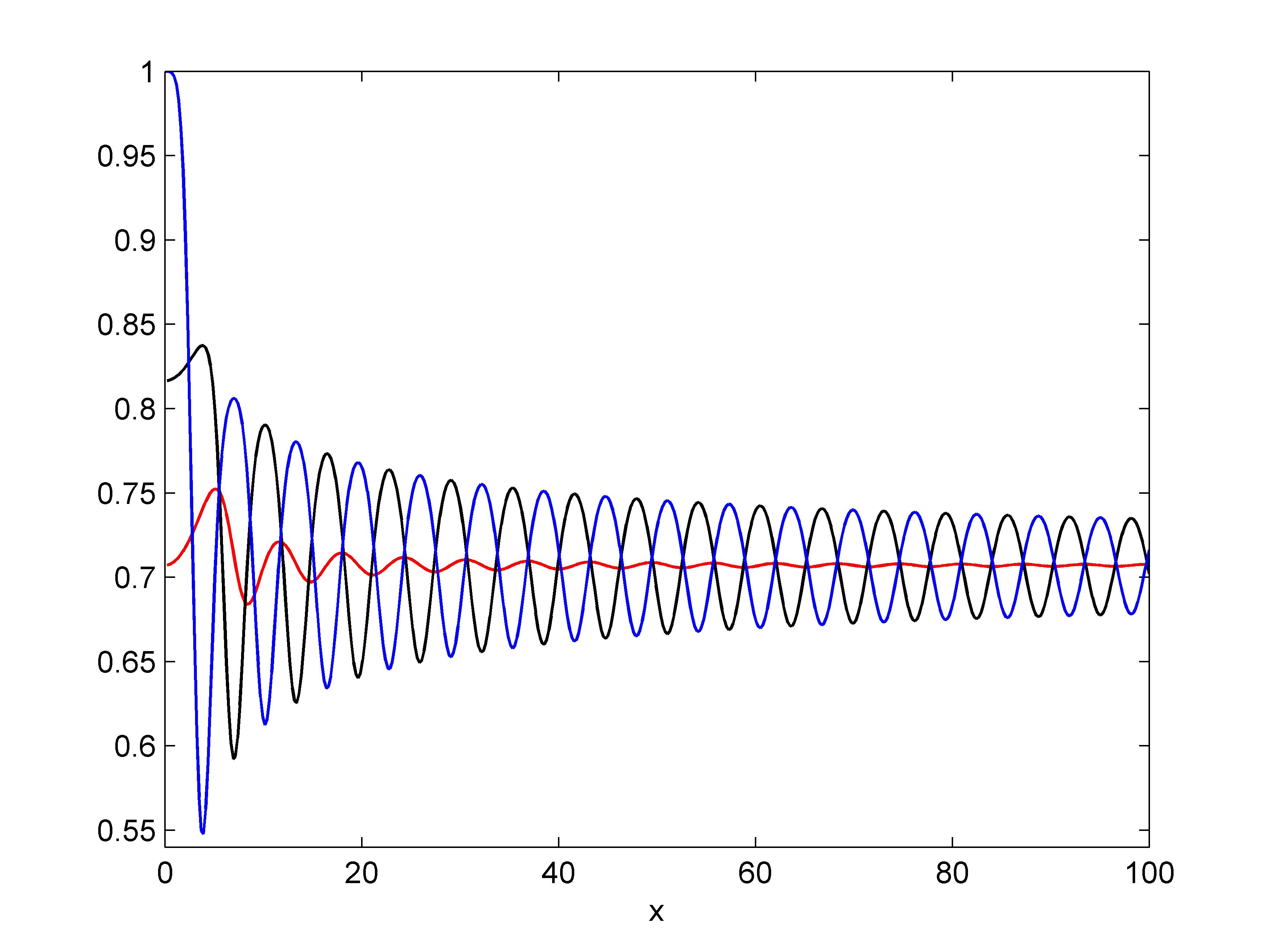}
\end{center}
\caption{Graphs of the functions \eqref{feq:1}--\eqref{feq:3}. See text for details.}
\label{fig}
\end{figure}

 Before moving on to the proof for the general case, it is interesting to actually compute the asymptotic lower bound to the critical radius for the two-dimensional case. The necessary information for this is contained in the three functions of in Figure \ref{fig}.


The red curve  in Figure \ref{fig}  is the plot of
\beq
f(x) \ \definedas \ \frac{1-\frac{2}{x}J_1(x)}{\sqrt{2-\frac{4}{x}J_1(x)-4((\frac{2}{x}J_1(x))')^2 }}.
\label{feq:1}
\eeq
The black one is the plot of
  \beq
g(x)\  \definedas \   \frac{1-J_0(x)}{\sqrt{2-2J_0(x)-2[J_0(x)']^2 }},
 \label{feq:2} \eeq
    and the blue one is the plot of
    \beq
 h(x)\  \definedas \    \frac{1+J_0(x)}{\sqrt{2+2J_0(x)-2[J_0(x)']^2 }}.
  \label{feq:3}  \eeq

However, from the proof of Theorem \ref{1} above we know that the
critical radius of the embeddings of  $i_n(S^2)$ in $\mathbb R^{(n+1)^2}$
is asymptotic to
$$\min\left\{\inf_{x\in [0,\infty)} f(x), \ \frac{1}{\sqrt{2}}\right\}.$$
 On the other hand, we know from (4.9)(4.10) in \cite{FA} that,  for $n$ odd where we have
 $\widetilde i_n(S^2)\cong S^2$  \eqref{af}, the limit is
  $$\min\left\{\inf_{x\in [0,\infty)}g(x), \   \inf_{x\in [0,\infty)} h(x), \ \frac{1}{\sqrt{2}}\right\};
  $$
  for $n$ even where the map is an immersion $\tilde i_n(S^2)\cong \mathbb RP^2$, the limit is
  $$\min\left\{  \inf_{x\in [0,\infty)} g(x), \ \frac{1}{\sqrt{2}}\right\}.
  $$

 It thus immediately follows from the Figure \ref{fig}  that the asymptotic lower bounds for the critical radii of the embeddings of $S^2$ by  $i_n$ we established here are larger than those for the embeddings $\widetilde i_n$ treated in \cite{FA}.
This is consistent with the expectations discussed in the Introduction.


\section{Proof of Theorem \ref{1} for the general case}\label{eeee}
The arguments for $S^2$ can be generalized to higher dimensions in a straightforward fashion, which we now describe.


The  kernel $K_n^d(x,y)$ of $\mathscr{H}_{\ell\leq n}^d$ can now be expressed as  \cite{bmo}
$$K_n^d(x,y)\ =\ \sum_{\ell=0}^n\sum_{j=1}^{k_\ell^d} \phi_j^\ell(x)\phi_j^\ell(y)\ =\ \frac{\pi_n^d/s_d}{{n+d/2 \choose n}}P_n^{(1+\lambda,\lambda)}(\cos\Theta(x,y)),$$
for $x,y\in S^d$,
where $\lambda=\frac{d-2}{2}$. 
 We first note that \cite{AH},
$$P_n^{(1+\lambda,\lambda)}(1) \ =\  {n+d/2 \choose n} \ =\ \frac{\Gamma(n+d/2+1)}{\Gamma(n+1)\Gamma(d/2+1)};$$and $$\ P_n^{(1+\lambda,\lambda)'}(1)\ = \frac{n+d}{2}{n+d/2 \choose n-1}.$$
We define the normalized spectral projection kernel as
$$\Pi_n^d(x,y)\ =\ \frac{s_d}{\pi^d_n}K_n^d(x,y)\ = \ {n+d/2 \choose n}^{-1}P_n^{(1+\lambda,\lambda)}(\Theta(x,y)). $$
Thus the norm of the map \eqref{id111} is $\|i_n^d(x)\|^2=\Pi_n^d(x,x)=1$, i.e.,
\begin{equation}
i_n^d:\,\,\, S^d\to  S^{\pi_n^d-1}.
\end{equation}
Following the arguments of \cite{FA}, the pull-back of the Euclidean metric is
\begin{equation}\label{metricback}(i_n^d)^*(g_E)\ =\ {n+d/2 \choose n}^{-1}P_n^{(1+\lambda,\lambda)'}(1)g_{S^d}\ =\ \frac{n(n+d)}{d+2} g_{S^d}.\end{equation}
Following the computations in \cite{FA}, the critical radius of the embedding is
\beqq
\inf_{\theta\in [0,\pi]} \frac{1-{n+d/2 \choose n}^{-1}P_n^{(1+\lambda,\lambda)}(\cos \theta)}{\sqrt{2-2{n+d/2 \choose n}^{-1}P_n^{(1+\lambda,\lambda)}(\cos\theta)-{n+d/2 \choose n}^{-1}\frac{ [P_n^{(1+\lambda,\lambda)'}(\cos\theta)\sin\theta]^2}{P_n^{(1+\lambda,\lambda)'}(1) }}}.\eeqq

We still have the following classical asymptotic estimate about the Jacobi polynomials.
Firstly, we have the asymptotic formula of Hilb's type below (\cite{Sz}, Theorem 8.21.12):
\begin{equation}\label{jacobi1}
\begin{split}
P_n^{(\alpha,\beta)}(\cos \theta)
&\ = \ \left(\sin \frac{\theta}{2}\right)^{-\alpha}\left(\cos \frac{\theta}{2}\right)^{-\beta}
\\ &  \qquad \times\left\{N^{-\alpha}\frac{\Gamma(n+\alpha+1)}{n!}
\left(\frac{\theta}{\sin \theta}\right)^{\frac{1}{2}}J_\alpha(N\theta)+R_{1,n}(\theta)\right\},
\end{split}
\end{equation}
where $N=n+(\alpha+\beta+1)/2$,
\begin{equation}
      R_{1,n}(\theta)\ = \
      \begin{cases}
        \theta^{\alpha+2}O(n^\alpha), &0\leq \theta\leq c/n,\\
        \theta^{\frac{1}{2}}O(n^{-\frac{3}{2}}),  &c/n\leq \theta\leq \pi-\epsilon,\\
      \end{cases}
\end{equation}
and $c,$ $\epsilon$ are uniform constants, independent of $n$. In our case $\alpha=1+\lambda$ and $\beta=\lambda$.\\
On the subinterval $c'/n\leq \theta \leq \pi-c'/n$, another asymptotic estimate is given by the Darboux formula (\cite{Sz}, Theorem 8.21.13)
\begin{equation}
P_n^{(1+\lambda,\lambda)}(\cos\theta)\ =\ \frac{1}{\sqrt{n}}k(\theta)\cos((n+\lambda+1)\theta+\gamma)+R_{2,n}(\theta),
\end{equation}
where $c$ is a  large enough, but uniform, constant, $\gamma=-(\lambda+{3}/{2}){\pi}/{2}$ and
\beqq
k(\theta)&=&\frac{1}{\sqrt{\pi}}\left(\sin\frac{\theta}{2}\right)^{-\lambda-3/2}\left(\cos\frac{\theta}{2}\right)^{-\lambda-1/2},
\\
R_{2,n}(\theta)&=&\frac{1}{\sqrt{n}}k(\theta)\frac{1}{n\sin\theta}O(1).
\eeqq

Near the end points, the asymptotic behavior of the Jacobi polynomials is given by the Mehler-Heine formulas (\cite{Sz}, p. 192)
\beqq
\lim_{n\to\infty}n^{-1-\lambda}P_n^{(1+\lambda,\lambda)}\left(\cos\frac{x}{n}\right)&=&\left(\frac{x}{2}\right)^{-1-\lambda}J_{1+\lambda}(x),
\\
\lim_{n\to\infty}n^{-\lambda}P_n^{(1+\lambda,\lambda)}\left(\cos\left(\pi-\frac{x}{n}\right)\right)&=&\left(\frac{x}{2}\right)^{-\lambda}J_{\lambda}(x),
\eeqq
where the limits are uniform on compact subsets of $\mathbb{R}$ and the $J_\nu$ are Bessel functions of the first kind with order $\nu$.

Again, following the arguments in Section \ref{dddd}, the global infimum is derived by breaking the argument up into the same four
subintervals used in the two-dimensional case and given in \eqref{eq:ints}.

Using the estimates given above, one can then argue as for the two-dimensional case. This leads to the following lower bound for the critical radius, as $n\to \infty$:
\beqq
\inf_{[0,\infty]}\frac{1-\left(\frac x2\right)^{-d/2}J_{d/2}(x)}
{\sqrt{2-\left(\frac{x}{2}\right)^{-d/2}J_{d/2}(x)
-2\Gamma(\frac{d}{2}+2)\left(\left(\left(\frac{x}{2}\right)^{-d/2}J_{d/2}(x)\right)'\right)^2 }}.
\eeqq
Again the expression does not diverge at $x=0$, as follows  by considering the full expansion of the Bessel function around $0$.
This completes the proof of Theorem \ref{1} for the general case. \qed

\end{document}